% Projet de Note aux CRAS, 3/3/2003
% 
% Benoit Fresse
% Fichier Plain-TeX, une compilation

% --- Definitions

% 1) Titles

% Fonts for titles
\font\sc=cmcsc10 \rm

% Numbers
\newcount\secnb
\newcount\thmnb
\newcount\parnb
\secnb=0
\parnb=0
\thmnb=64

% Vertical skips
\def\smallskip{\par\vskip 2mm}
\def\medskip{\par\vskip 5mm}
\def\goodbreak{\penalty -100}

% Section
\def\section#1{\global\advance\secnb by 1
\parnb=0
\medskip\goodbreak
\centerline{\bf\S\the\secnb.\ #1}
\smallskip\nobreak}

% Bibliography
\def\references{\medskip\goodbreak
\centerline{\bf R\'ef\'erences}
\smallskip\nobreak}

% Paragraph
\def\paragraph#1{\global\advance\parnb by 1
\smallskip\goodbreak
\noindent{\bf\the\parnb)\ #1}
\par\nobreak}

% Theorem
\def\thm#1{\global\advance\thmnb by 1
\smallskip\goodbreak
\noindent{\sc #1}\ \char\the\thmnb
\par\nobreak}

% 2) Cross references

\newtoks\refid

% Label
\def\label#1{\relax}

% Reference
\def\ref#1{\refid={\csname crossref#1\endcsname}
\the\refid}

% 3) Bibliography

% Reference number
\newcount\bibrefnb
\bibrefnb=0

% Label of a reference
\def\bibitem#1{\global\advance\bibrefnb by 1\item{\the\bibrefnb.}}
\def\authorname#1{{\sc #1}}

% Reference to the bibliography
\def\cite#1{\refid={\csname bibref#1\endcsname}
[\the\refid]}

% 4) Last page

\def\bye{\medskip\noindent{\it Addresse:} Benoit Fresse\par
Laboratoire J.A. Dieudonn\'e\par
Universit\'e de Nice-Sophia-Antipolis\par
Parc Valrose\par
F-06108 Nice Cedex 02 (France)\par
\noindent{\it Courriel:} {\tt fresse@math.unice.fr}\par
\vfill\end}

% 5) Abstracts

% Abstracts
\def\beginresume{\smallskip
\begingroup\leftskip=1cm\rightskip=1cm
{\sc R\'esum\'e}. --- }
\def\endresume{\endgroup}

\def\beginabstract{\smallskip
\begingroup\leftskip=1cm\rightskip=1cm
{\sc Abstract}. --- }
\def\endabstract{\endgroup}

% 5) Notations

\def\N{{\bf N}}
\def\F{{\bf F}}

\def\E{{\cal E}}
\def\Free{{\cal F}}

\def\X{{\cal X}}
\def\C{{\cal C}}

\def\TR{\mathop{TR}\nolimits}
\def\id{\mathop{id}\nolimits}
\def\dg{\mathop{dg}\nolimits}

\def\LMod#1{\mathop{#1{\rm Mod}}\nolimits}
\def\LAlg#1{\mathop{#1{\rm Alg}}\nolimits}

% --- Text

\hsize=15cm
\vsize=20.5cm

\headline={\ifnum\pageno=1{\it Projet de note aux C.R.Acad.Sci.Paris S\'er. I Math. (Topologie)}\hfil\else\hfil\fi}

% Titre
\vbox{\medskip
\centerline{\bf La construction bar d'une alg\`ebre comme alg\`ebre de Hopf E-infini}
\smallskip
\centerline{\sc Benoit Fresse}
\smallskip
\centerline{3/3/2003}
\medskip
\beginresume
On prouve que la construction bar d'une alg\`ebre $E_\infty$ forme une alg\`ebre $E_\infty$.
Plus pr\'ecis\'ement,
on montre que la construction bar d'une alg\`ebre sur l'op\'erade des surjections
poss\`ede une structure d'alg\`ebre de Hopf sur l'op\'erade de Barratt-Eccles.
(L'op\'erade des surjections et l'op\'erade de Barratt-Eccles sont des op\'erades $E_\infty$ classiques.)
\endresume
\medskip
\centerline{\bf The bar construction of an algebra as an E-infinite Hopf algebra}
\smallskip
\beginabstract
We prove that the bar construction of an $E_\infty$ algebra forms an $E_\infty$ algebra.
To be more precise,
we provide the bar construction of an algebra over the surjection operad
with the structure of a Hopf algebra over the Barratt-Eccles operad.
(The surjection operad and the Barratt-Eccles operad are classical $E_\infty$ operads.)
\endabstract
\smallskip}

\medskip
\noindent{\bf Abridged English Version}
\medskip

We fix a ground field $\F$ of characteristic $2$.
We let $\Sigma_r$, $r\in\N$, denote the sequence of symmetric groups.
We consider operads in the category $\dg\LMod{\F}$ of differential graded modules over $\F$
(for short {\it dg-modules}\/).

We denote the operad of associative and commutative algebras by the letter $\C$.
We recall that an $E_\infty$ operad consists of a dg-operad $\Free$ quasi-isomorphic to $\C$
and whose components $\Free(r)$, $r\in\N$, are projective complexes of $\Sigma_r$-modules.
The category of algebras over a fixed $E_\infty$ operad $\LAlg{\Free}$ is equipped with the structure
of a semi-model category ({\it cf}. \cite{Mandell}).
The purpose of this note is to make explicit a model of the suspension of an algebra in $\LAlg{\Free}$.

This model is given by the classical bar construction of associative algebras $\bar{B}(A)$.
More specifically, the bar construction of an associative and commutative algebra
is equipped with the {\it shuffle} product and forms an associative and commutative algebra.
We extend this construction to the context of $E_\infty$ algebras
by introducing particular $E_\infty$ operads.
Namely: the {\it surjection operad} $\X$
and the {\it Barratt-Eccles operad} $\E$.

The components $\X(2)$ and $\E(2)$ of these operads are both isomorphic to the classical free resolution
of the trivial representation of $\Sigma_2$.
To be more explicit,
we have homogeneous elements $\theta_d$ such that $\E(2)_d = \X(2)_d = \F[\Sigma_2]\cdot\theta_d$.
In addition,
the differential of the dg-modules $\E(2)_* = \X(2)_*$
verifies the formula $\delta(\theta_d) = \theta_{d-1} + \tau\cdot\theta_{d-1}$,
where $\tau$ denotes the transposition of $\Sigma_2$.
For a given $\X$-algebra $A$,
the operation $\theta_d: A\otimes A\,\rightarrow\,A$ associated to $\theta_d\in\X(2)$
is also denoted by $a_1\smile_d a_2 = \theta_d(a_1,a_2)$.
In fact,
the element $\theta_0\in\X(2)_0$
satisfies the relation
$\theta_0(\theta_0,1) = \theta_0(1,\theta_0)$
in the surjection operad.
Accordingly,
the product $\smile_0$ is associative in $A$.
Similarly,
the product $\smile_d$ gives rise to the boundary relation
$a_1\smile_{d-1} a_2 + a_2\smile_{d-1} a_1 = \delta(a_1\smile_d a_2) + \delta(a_1)\smile_d a_2 + a_1\smile_d\delta(a_2)$,
because we have $\delta(\theta_d) = \theta_{d-1} + \tau\cdot\theta_{d-1}$
in the surjection operad.

The work of Baues ({\it cf}. \cite{Baues}) proves that the bar construction of an $\X$-algebra
is equipped with an associative product.
There is also a sequence of products $\smile_d: \bar{B}(A)\otimes\bar{B}(A)\,\rightarrow\,\bar{B}(A)$
defined by Kadeishvili in the article \cite{KadeishviliHigherProducts}
and that verify the boundary relation above.
We generalize Kadeishvili's construction
and we obtain the following theorem:

\thm{Theorem}
{\it Let $A$ be an algebra over the surjection operad $\X$.
The bar construction $\bar{B}(A)$ is equipped with the structure of a Hopf algebra over the Barratt-Eccles operad $\E$
such that the operation $\theta_d: \bar{B}(A)\otimes\bar{B}(A)\,\rightarrow\,\bar{B}(A)$
associated to the element $\theta_d\in\E(2)_d$
agrees with Kadeishvili's product $\smile_d$.
The bar construction $\bar{B}(A)$ together with this structure forms a cogroup object in the homotopy category of $\E$-algebras
and is equivalent to the suspension of $A$.}

\smallskip
We define an operad morphism $\TR: \E\,\rightarrow\,\X$ in \cite{BergerFressePrisms} and \cite{BergerFresseCochains}.
This morphism gives any algebra over $\X$ the structure of an algebra over $\E$.
Therefore, it makes sense to consider the suspension of an $\X$-algebra
in the category of $\E$-algebras.

The normalized cochain complex of a simplicial set $A = N^*(X)$ is equipped with the structure of an algebra
over the surjection operad $\X$
({\it cf}. \cite{BergerFresseCochains}, \cite{McClureSmithSurjections}).
Moreover, according to a general result of Mandell ({\it cf}. \cite{Mandell}),
this structure determines the $2$-adic homotopy type of $X$
(since we consider cochains with $\F = \F_2$ coefficients).
One proves that the suspension of $N^*(X)$
in the homotopy category of $\E$-algebras is equivalent to $N^*(\Omega X)$,
the cochain algebra of the loop space of $X$
({\it cf}. \cite{FresseMaps}, \cite{Mandell}).
Consequently:

\thm{Theorem}
{\it We assume that $X$ is a pointed connected simplicial set such that $\pi_1(X)$ is a finite $2$-group
and $H^*(X,\F_2)$ is a finite dimensional $\F_2$-module for all $*\in\N$.
If $A = N^*(X)$, the cochain algebra of $X$,
then, in the homotopy category of $\E$-algebras,
the bar construction $\bar{B}(N^*(X))$ is equivalent to $N^*(\Omega X)$,
the cochain algebra of the loop space of $X$.}

\smallskip
The works of Smirnov ({\it cf}. \cite{SmirnovCLoop}), Justin R. Smith ({\it cf}. \cite{JustinSmithOperads})
and Kadeishvili-Saneblidze ({\it cf}. \cite{KadeishviliSaneblidze})
predict the existence of an $E_\infty$ structure on the bar construction $\bar{B}(N^*(X))$.
Our theorems make this structure explicit.

\medskip
\hbox to 3cm{\hrulefill}
\thmnb=64

\medskip
On fixe un corps de base $\F$ de caract\'eristique $2$.
On consid\`ere des op\'erades dans la cat\'egorie $\dg\LMod{\F}$
des modules diff\'erentiels gradu\'es sur $\F$
(en abr\'eg\'es {\it dg-modules}\/).
On note $\Sigma_r$, $r\in\N$, la suite des groupes de permutations.

\section{R\'esultats}

\smallskip
L'op\'erade associ\'ee aux alg\`ebres associatives et commutatives
est d\'esign\'ee par la lettre $\C$.
On rappelle qu'une {\it op\'erade $E_\infty$} est une dg-op\'erade $\Free$ quasi-isomorphe \`a $\C$
et dont les composantes $\Free(r)$, $r\in\N$,
forment des complexes projectifs de $\Sigma_r$-modules.
La cat\'egorie d'alg\`ebres associ\'ee \`a une telle op\'erade $\LAlg{\Free}$
poss\`ede une semi-structure mod\`ele naturelle ({\it cf}. \cite{Mandell}).
Le but de cette note est de donner un mod\`ele explicite de la suspension d'une alg\`ebre dans $\LAlg{\Free}$.

Ce mod\`ele est fourni par la construction bar classique des alg\`ebres associatives $\bar{B}(A)$.
Plus sp\'ecifiquement, on sait que le produit {\it shuffle}\/ donne \`a la construction bar d'une alg\`ebre commutative
la structure d'une alg\`ebre commutative.
On \'etend cette construction au cadre des alg\`ebres $E_\infty$
en introduisant des op\'erades $E_\infty$ particuli\`eres
qui sont l'{\it op\'erade des surjections} $\X$
et l'{\it op\'erade de Barratt-Eccles} $\E$.

Les composantes $\X(2)$ et $\E(2)$ de ces op\'erades sont isomorphes \`a la r\'esolution libre classique
de la repr\'esentation triviale de $\Sigma_2$.
Plus explicitement, on a une suite d'\'el\'ements homog\`enes $\theta_d$
tels que $\E(2)_d = \X(2)_d = \F[\Sigma_2]\cdot\theta_d$.
En outre,
la diff\'erentielle des dg-modules $\E(2)_* = \X(2)_*$
v\'erifie la formule $\delta(\theta_d) = \theta_{d-1} + \tau\cdot\theta_{d-1}$,
en notant $\tau$ la transposition de $\Sigma_2$.
Si $A$ est une $\X$-alg\`ebre,
alors l'op\'eration $\theta_d: A\otimes A\,\rightarrow\,A$ associ\'ee \`a l'\'el\'ement $\theta_d\in\X(2)_d$
est \'egalement not\'ee $a_1\smile_d a_2 = \theta_d(a_1,a_2)$.
En fait,
l'\'el\'ement $\theta_0$ satisfait la relation $\theta_0(\theta_0,1) = \theta_0(1,\theta_0)$
dans l'op\'erade des surjections.
Le produit correspondant $\smile_0$ est donc associatif sur $A$.
De m\^eme,
comme $\theta_d$ a pour diff\'erentielle $\delta(\theta_d) = \theta_{d-1} + \tau\cdot\theta_{d-1}$,
le produit sup\'erieur $\smile_d$ donne lieu \`a la relation de bord
$a_1\smile_{d-1} a_2 + a_2\smile_{d-1} a_1 = \delta(a_1\smile_d a_2) + \delta(a_1)\smile_d a_2 + a_1\smile_d\delta(a_2)$.

Les travaux de Baues ({\it cf}. \cite{Baues}) montrent que la construction bar d'une $\X$-alg\`ebre $\bar{B}(A)$
poss\`ede une structure d'alg\`ebre associative.
On a aussi une suite de produits
$\smile_d: \bar{B}(A)\otimes\bar{B}(A)\,\rightarrow\,\bar{B}(A)$
(v\'erifiant la relation de bord ci-dessus)
que Kadeishvili d\'efinit de fa\c con explicite dans l'article \cite{KadeishviliHigherProducts}.
On \'etend la construction de Kadeishvili pour associer une op\'eration sur la construction bar
\`a tout \'el\'ement de l'op\'erade de Barratt-Eccles $\E$.
On obtient ainsi le th\'eor\`eme\break suivant :

\thm{Th\'eor\`eme}\label{MainTheorem}
{\it Si $A$ est une alg\`ebre sur l'op\'erade des surjections $\X$,
alors la construction bar $\bar{B}(A)$ poss\`ede une structure naturelle d'alg\`ebre
de Hopf sur l'op\'erade de Barratt-Eccles $\E$
telle que l'op\'eration
$\theta_d: \bar{B}(A)\otimes\bar{B}(A)\,\rightarrow\,\bar{B}(A)$
associ\'ee \`a l'\'el\'ement
$\theta_d\in\E(2)_d$
est le produit $\smile_d$ d\'efini par Kadeishvili.
Quand elle est munie de cette structure,
la construction bar $\bar{B}(A)$ d\'efinit un objet en cogroupe dans la cat\'egorie homotopique des $\E$-alg\`ebres
et est \'equivalente \`a la suspension de $A$.}

\smallskip
On d\'efinit un morphisme d'op\'erades $\TR: \E\,\rightarrow\,\X$ dans les articles \cite{BergerFressePrisms} et \cite{BergerFresseCochains}.
Ce morphisme fait de toute alg\`ebre sur $\X$ une alg\`ebre sur $\E$ par restriction de structure.
C'est pourquoi on peut parler de la suspension d'une $\X$-alg\`ebre dans la cat\'egorie des $\E$-alg\`ebres.

On sait que le complexe des cocha\^\i nes normalis\'ees d'un ensemble simplicial $A = N^*(X)$
poss\`ede une structure naturelle d'alg\`ebre sur l'op\'erade des surjections $\X$
({\it cf}. \cite{BergerFresseCochains}, \cite{McClureSmithSurjections}).
De plus, d'apr\`es un r\'esultat g\'en\'eral de Mandell ({\it cf}. \cite{Mandell}),
cette structure suffit \`a d\'eterminer le type d'homotopie $2$-adique de $X$
(quand on prend $\F = \F_2$ comme corps de coefficients).
On montre que la suspension de $N^*(X)$
dans la cat\'egorie homotopique des $\E$-alg\`ebres
est \'equivalente \`a $N^*(\Omega X)$,
l'alg\`ebre des cocha\^\i nes de l'espace des lacets de $X$
({\it cf}. \cite{FresseMaps}, \cite{Mandell}).
En cons\'equence :

\thm{Th\'eor\`eme}
{\it On suppose que $X$ est un ensemble simplicial point\'e connexe dont la cohomologie $H^*(X,\F_2)$ est finie en tout degr\'e
et tel que $\pi_1(X)$ forme un $2$-groupe fini.
Si $A = N^*(X)$, l'alg\`ebre des cocha\^\i nes de $X$,
alors la construction bar $\bar{B}(N^*(X))$ est \'equivalente
dans la cat\'egorie homotopique des $\E$-alg\`ebres
\`a $N^*(\Omega X)$, l'alg\`ebre des cocha\^\i nes de l'espace des lacets de $X$.}

\smallskip
L'existence d'une structure $E_\infty$ sur la construction bar $\bar{B}(N^*(X))$
est assur\'ee par les travaux de Smirnov ({\it cf}. \cite{SmirnovCLoop}),
de Justin R. Smith ({\it cf}. \cite{JustinSmithOperads})
et de Kadeishvili-Saneblidze ({\it cf}. \cite{KadeishviliSaneblidze}).
Nos th\'eor\`emes rendent une telle structure explicite.
Le but de la seconde partie de cette note est de d\'efinir l'op\'eration $w: \bar{B}(A)^{\otimes r}\,\rightarrow\,\bar{B}(A)$
associ\'ee \`a un \'el\'ement $w\in\E(r)$.
La d\'emonstration du th\'eor\`eme \ref{MainTheorem} sera publi\'ee ult\'erieurement.

\section{Construction des op\'erations sur la construction bar}

\smallskip
On reprend les conventions classiques classique du calcul diff\'erentiel gradu\'e.
Un dg-module $V$ est gradu\'e inf\'erieurement $V = V_*$ ou sup\'erieurement $V = V^*$,
la relation $V_{d} = V^{-d}$
rendant une graduation inf\'erieure \'equivalente \`a une graduation sup\'erieure.
La diff\'erentielle d'un dg-module est g\'en\'eralement not\'ee
$\delta: V_*\,\rightarrow\,V_{*-1}$.

\paragraph{Rappels sur l'op\'erade de Barratt-Eccles et l'op\'erade des surjections}
On reprend les conventions des articles \cite{BergerFressePrisms} et \cite{BergerFresseCochains}.
On rappelle que l'op\'erade de Barratt-Eccles $\E$ est d\'efinie
par la construction bar homog\`ene normalis\'ee
des groupes sym\'etriques $\Sigma_r$.
Explicitement,
le module $\E(r)$ est engendr\'e en degr\'e $d$
par les $d+1$-uplets non-d\'eg\'en\'er\'es de permutations
$(w_0,\ldots,w_d)\in\Sigma_r\times\cdots\times\Sigma_r$.
Un simplexe $w = (w_0,\ldots,w_d)$ est d\'eg\'en\'er\'e (et repr\'esente $0$ dans $\E(r)_d$)
si on a $w_{i+1}=w_{i}$ pour quelque $i\in\{0,\ldots,d-1\}$.
Ainsi,
l'\'el\'ement $\theta_d\in\E(2)_d$ correspondant au produit $\smile_d$ est repr\'esent\'e
par le simplexe alternant $\theta_d = (\id,\tau,\id,\tau,\ldots)$.
La diff\'erentielle de $\E(r)$ est d\'efinie par la formule classique
$\delta(w_0,\ldots,w_d) = \sum_{i=0}^d (w_0,\ldots,\widehat{w_i},\ldots,w_d)$.

Les composantes $\X(r)_d$ de l'op\'erade des surjections $\X$ sont engendr\'ees par les surjections non-d\'eg\'en\'er\'ees
$u: \{1,\ldots,r+d\}\,\rightarrow\,\{1,\ldots,r\}$,
une surjection $u$ \'etant d\'eg\'en\'er\'ee si $u(i+1)=u(i)$ pour quelque $i\in\{1,\ldots,r+d-1\}$
(auquel cas, on suppose que $u$ repr\'esente $0$ dans $\X(r)_d$).
Une surjection $u\in\X(r)_d$ est d\'etermin\'ee par la suite de ses valeurs $\underline{u} =  (u(1),\ldots,u(r+d))$.

On d\'efinit dans l'article \cite{BergerFresseCochains}
une certaine d\'ecomposition de $\underline{u}$
en sous-suites $\underline{u}_0,\ldots,\underline{u}_d$
(ce sont les {\it lignes}\/ de l'{\it arrangement en table}\/ de $u\in\X(r)_d$).
Pour caract\'eriser cette d\'ecomposition, on sp\'ecifie les termes de $\underline{u}$
qui d\'efinissent les derniers \'el\'ements de $\underline{u}_0,\ldots,\underline{u}_{d-1}$
(les {\it c\'esures}\/ de la surjection $u$) :
ce sont les termes de $\underline{u}$
qui ne forment pas la derni\`ere occurrence d'une valeur $k = 1,\ldots,r$
dans la suite $\underline{u}$.
Par exemple, pour $\underline{u} = (1,4,\underline{2},5,\underline{3},2,3)$, on obtient :
$$\underline{u} = (\underbrace{1,4,2}_{\underline{u}_0};\underbrace{5,3}_{\underline{u}_1};\underbrace{2,3}_{\underline{u}_2})$$
(les c\'esures sont soulign\'ees dans la suite initiale).

\paragraph{Rappels sur la construction bar}
La {\it cog\`ebre tensorielle}\/ engendr\'ee par un dg-module $V$, not\'ee $T^c(V)$,
est form\'ee par le module $T^c(V) = \bigoplus_{n=0}^{\infty} V^{\otimes n}$
muni de la diagonale
$\Delta: T^c(V)\,\rightarrow\,T^c(V)\otimes T^c(V)$
d\'efinie par la d\'econcat\'enation des tenseurs.
On obtient ainsi une structure de cog\`ebre associative.

On suppose que $A$ est une alg\`ebre augment\'ee sur l'op\'erade des surjections $\X$.
On note $\tilde{A}$ l'id\'eal d'augmentation de $A$.
On consid\`ere la suspension de $\tilde{A}$ (dans la cat\'egorie des dg-modules)
dont les composantes homog\`enes sont d\'efinies par la relation
$(\Sigma\tilde{A})^* = \tilde{A}^{*+1}$.
La construction bar $\bar{B}(A)$ est d\'efinie par la cog\`ebre tensorielle $\bar{B}(A) = T^c(\Sigma\tilde{A})$
munie d'une diff\'erentielle $b': \bar{B}(A)\,\rightarrow\,\bar{B}(A)$
qui est d\'etermin\'ee par le produit associatif de $A$
associ\'e \`a l'op\'eration $\theta_0\in\X(2)$.
Explicitement,
cette diff\'erentielle $b': \bar{B}(A)\,\rightarrow\,\bar{B}(A)$
est donn\'ee par la formule
$$b'(\Sigma a_1\otimes\cdots\otimes\Sigma a_n)
= \sum_{i=1}^{n-1} \Sigma a_1\otimes\cdots\otimes\Sigma(a_i\smile_0 a_{i+1})\otimes\cdots\otimes\Sigma a_n.$$

Dans le prochain paragraphe,
on associe \`a tout \'el\'ement $w\in\E(r)$ une application $\tilde{w}: T^c(\Sigma\tilde{A})^{\otimes r}\,\rightarrow\,\tilde{A}$.
On \'etend cette application
en une op\'eration $w: T^c(\Sigma\tilde{A})^{\otimes r}\,\rightarrow\,T(\Sigma\tilde{A})$
(sur la construction bar)
en utilisant la structure de cog\`ebre de $T^c(\Sigma\tilde{A})$.
On rappelle que l'op\'erade de Barratt-Eccles est munie d'une diagonale coassociative $\Delta: \E(r)\,\rightarrow\,\E(r)\otimes\E(r)$ comme la cog\`ebre tensorielle.
(En cons\'equence, les modules $\E(r)$ forment une op\'erade dans la cat\'egorie des cog\`ebres ;
on dit aussi que l'op\'erade de Barratt-Eccles $\E$ est une {\it op\'erade de Hopf}\/).
On note $\Delta^n(w) = \sum\nolimits_i w^i_{(1)}\otimes\cdots\otimes w^i_{(n)}$ la diagonale it\'er\'ee
de l'\'el\'ement $w\in\E(r)$ dans $\E(r)^{\otimes n}$.
On pose explicitement
$$w(c_1,\ldots,c_r) = \sum\nolimits_i \Sigma\widetilde{w}^i_{(1)}(c^i_{1(1)},\ldots,c^i_{r(1)})\otimes\cdots
\otimes\Sigma\widetilde{w}^i_{(n)}(c^i_{1(n)},\ldots,c^i_{r(n)}),$$
en notant
$$\Delta^n(c_k) = \sum\nolimits_i c^i_{k(1)}\otimes\cdots\otimes c^i_{k(n)}\in T^c(\Sigma\tilde{A})^{\otimes n},$$
les diagonales it\'er\'ees des tenseurs $c_1,\ldots,c_r\in T^c(\Sigma\tilde{A})$.
On montre que cette construction donne \`a la construction bar $\bar{B}(A) = T^c(\Sigma\tilde{A})$ une structure de $\E$-alg\`ebre.
(On obtient en fait une $\E$-alg\`ebre dans la cat\'egorie des dg-cog\`ebres ;
c'est pourquoi on dit que $\bar{B}(A)$ forme une {\it alg\`ebre de Hopf}\/ sur $\E$).

\paragraph{Construction des op\'erations sur la construction bar}
On d\'efinit l'application
$\widetilde{w}: T^c(\Sigma\tilde{A})\otimes\cdots\otimes T^c(\Sigma\tilde{A})\,\rightarrow\,\tilde{A}$
associ\'ee \`a un simplexe $w = (w_0,\ldots,w_d)\in\E(r)_d$.
Cette application est donn\'ee sur chaque composante
$(\Sigma\tilde{A})^{\otimes n_1}\otimes\cdots\otimes(\Sigma\tilde{A})^{\otimes n_r}
\subset T^c(\Sigma\tilde{A})\otimes\cdots\otimes T^c(\Sigma\tilde{A})$
par une somme d'op\'erations
$u: A^{\otimes n_1}\otimes\cdots\otimes A^{\otimes n_r}\,\rightarrow\,A$
associ\'ees \`a des surjections $u\in\X(n_1+\cdots+n_r)$
{\it admissibles} par rapport \`a $w$.

Dans le contexte de cette construction,
il est naturel de repr\'esenter une surjection $u\in\X(n_1+\cdots+n_r)$ par une application \`a valeurs
dans l'ensemble $M = \{1_1,2_1,\ldots,(n_1)_1,\ldots,1_r,2_r,\ldots,(n_r)_r\}$
constitu\'e de $r$ intervalles d'entiers.
On note que chaque permutation $w_i$ de $w = (w_0,\ldots,w_d)$
d\'efinit un ordre sur les intervalles de l'ensemble $M$.
Explicitement,
si on se donne $k_s,l_t\in M$ avec $s\not=t$,
alors on \'ecrit $k_s<_{w_i} l_t$
quand le couple $(s,t)$ constitue une sous-suite de $\underline{w}_i = (w_i(1),\ldots,w_i(r))$.
On a par exemple
$$k_1<_{(1,3,2)} m_3<_{(1,3,2)} l_2,$$
quelque soient $k_1\in\{1_1,2_1,\ldots,(n_1)_1\}$, $l_2\in\{1_2,2_2,\ldots,(n_2)_2\}$
et $m_3\in\{1_3,2_3,\ldots,(n_3)_3\}$.

Une surjection $u$ est admissible par rapport \`a $w$
quand on peut r\'epartir les lignes de l'arran\-ge\-ment en table de $u$
en $d+1$ blocs
$$\underline{u} = \,\underbrace{\underline{u}^0_1,\ldots,\underline{u}^0_{e_0}}_{\underline{u}^0}\,,
\,\underbrace{\underline{u}^1_0,\underline{u}^1_1,\ldots,\underline{u}^1_{e_1}}_{\underline{u}^1}\,,
\,\ldots\,,
\,\underbrace{\underline{u}^d_0,\underline{u}^d_1,\ldots,\underline{u}^d_{e_d}}_{\underline{u}^d}$$
tout en respectant les propri\'et\'es caract\'eristiques suivantes.
On consid\`ere la ligne $\underline{u}^i_j$ de cet arrangement.
On note $(k_1)_{s_1},\ldots,(k_m)_{s_m}$ les valeurs des c\'esures des lignes pr\'ec\'edant $\underline{u}^i_j$
dont l'occurrence finale n'appara\^\i t pas d\'ej\`a dans la table
(au dessus de la ligne $\underline{u}^i_j$).
On suppose que ces valeurs sont ordonn\'ees selon l'ordre d\'efini par la permutation $w_i$.
On a explicitement $(k_1)_{s_1}\,<_{w_i}\,\cdots\,<_{w_i}\,(k_m)_{s_m}$.
(On observera que, par construction, un nombre donn\'e $s\in\{1,\ldots,r\}$ appara\^\i t toujours au plus une fois dans la suite $s_1,\ldots,s_m$.)
Si $j = 0$ (on suppose donc que $\underline{u}^i_j$ est la premi\`ere ligne du bloc $\underline{u}^i$),
alors on demande que la ligne $\underline{u}^i_j$
soit constitu\'ee par la suite $\underline{u}^i_j = ((k_1)_{s_1},\ldots,(k_{l})_{s_{l}})$,
avec $1\leq l\leq m$.
Sinon (si $j > 0$),
on demande \`a avoir $\underline{u}^i_j = ((k_0)_{s_0},(k_1)_{s_1},\ldots,(k_{l})_{s_{l}})$,
le premier terme de cette ligne repr\'esentant la premi\`ere occurrence d'une valeur $(k_0)_{s_0}\in M$ dans la suite $\underline{u}$
et v\'erifiant $(k_0)_{s_0}<_{w_i}\,(k_1)_{s_1}$.
On demande aussi que les valeurs $k_s\in M$ associ\'ees \`a un nombre $s\in\{1,\ldots,r\}$ fix\'e
apparaissent dans l'ordre croissant $k = 1,\ldots,n_s$
dans la suite $\underline{u}$.

On note bien que le premier \'el\'ement d'une ligne $\underline{u}^i_j$ telle que $j = 1,\ldots,e_i$
repr\'esente la premi\`ere occurrence d'une valeur $k_s\in M$ dans $\underline{u}$.
La propri\'et\'e inverse caract\'erise donc la premi\`ere ligne d'un bloc de notre d\'ecomposition.

\paragraph{Exemples}
On peut d\'eterminer facilement les surjections $u\in\X(p+q)$ qui sont admissibles par rapport aux op\'erations $\theta_d\in\E(2)_d$.
Ainsi, pour $\theta_0 = (\id)\in\E(2)_0$, les surjections admissibles sont de la forme
$$\underline{u} = (\,\underbrace{1_2}_{\underline{u}^0_1}\,,\,\underbrace{1_1,1_2}_{\underline{u}^0_2}\,,
\,\underbrace{2_1,1_2}_{\underline{u}^0_3}\,,\,\underbrace{3_1,1_2}_{\underline{u}^0_4}\,,
\,\ldots\,,\,\underbrace{p_1,1_2}_{\underline{u}^0_{p+1}}\,).$$
On reconnait les \'el\'ements de $\X(p+1)$ associ\'es aux op\'erations ``{\it braces}'' de Getzler-Kadeishvili
({\it cf}. \cite{BergerFresseCochains}, \cite{McClureSmithPrisms}, \cite{McClureSmithSurjections}).
Pour $\theta_1 = (\id,\tau)\in\E(2)_1$, on obtient les surjections de la forme
$$(\,\underbrace{1_2}_{\underline{u}^0_1}\,,\,\underbrace{1_1,1_2}_{\underline{u}^0_2}\,,
\,\underbrace{2_1,1_2}_{\underline{u}^0_3}\,,\,\underbrace{3_1,1_2}_{\underline{u}^0_4}\,,
\,\ldots\,,\,\underbrace{p_1}_{\underline{u}^0_{p+1}}\,,
\,\underbrace{1_2,p_1}_{\underline{u}^1_0}\,,
\,\underbrace{2_2,p_1}_{\underline{u}^1_1}\,,\,\underbrace{3_2,p_1}_{\underline{u}^1_2}\,,
\,\ldots\,,\,\underbrace{q_2,p_1}_{\underline{u}^1_{q}}\,).$$
On reconnait les \'el\'ements $E^1_{p q}$ introduits par Kadeishvili dans l'article \cite{KadeishviliHigherProducts}.

\references{\parindent=0cm\leftskip=1cm\rightskip=1cm

\bibitem{Baues}\authorname{H. Baues},
Geometry of loop spaces and the cobar construction,
Memoirs of the American Mathematical Society {\bf 25},
American Mathematical Society, 1980.

\bibitem{BergerFressePrisms}\authorname{C. Berger, B. Fresse},
{\it Une d\'ecomposition prismatique de l'op\'erade de Barratt-Eccles},
C. R. Acad. Sci. Paris S\'er. I {\bf 335} (2002), 365-370.

\bibitem{BergerFresseCochains}\authorname{C. Berger, B. Fresse},
{\it Combinatorial operad actions on cochains},
pr\'epublication {\tt arXiv:math.AT/0109158} (2001).

\bibitem{FresseMaps}\authorname{B. Fresse},
{\it Derived division functors and mapping spaces},
pr\'epublication {\tt arXiv:\break math.AT/0208091} (2001).

\bibitem{KadeishviliHigherProducts}\authorname{T. Kadeishvili},
{\it Cochain operations defining Steenrod $\smile_i$-products in the bar construction},
pr\'epublication {\tt arXiv:math.AT/0207010} (2002).

\bibitem{KadeishviliSaneblidze}\authorname{T. Kadeishvili, S. Saneblidze},
{\it Iterating the bar construction},
Georgian Math. J. {\bf 5} (1998), 441-452.

\bibitem{Mandell}\authorname{M. Mandell},
{\it $E_\infty$ algebras and $p$-adic homotopy theory},
Topology {\bf 40} (2001), 43-94.

\bibitem{McClureSmithPrisms}\authorname{J. McClure, J.H. Smith},
{\it A solution of Deligne's Hochschild cohomology conjecture}
{\it in} ``Recent progress in homotopy theory (Baltimore, 2000)'',
Contemp. Math. {\bf 293}, Amer. Math. Soc. (2002), 153-193.

\bibitem{McClureSmithSurjections}\authorname{J. McClure, J.H. Smith},
{\it Multivariable cochain operations and little $n$-cubes},
pr\'e\-pu\-bli\-ca\-tion {\tt arXiv:math.QA/0106024} (2001).

\bibitem{SmirnovCLoop}\authorname{V. Smirnov},
{\it On the chain complex of an iterated loop space},
Izv. Akad. Nauk SSSR Ser. Mat. {\bf 53} (1989), 1108-1119.
English translation in Math. USSR-Izv. {\bf 35} (1990), 445-455.

\bibitem{JustinSmithOperads}\authorname{J.R. Smith},
{\it Operads and algebraic homotopy},
pr\'epublication {\tt arXiv:math.AT/\break 0004003} (2000).

}

\bye